\documentclass[5p,twocolumn,sort&compress,times,preprint,normalem]{elsarticle}

\usepackage{graphics}
\usepackage{graphicx}

\usepackage{amsfonts, amsthm, amssymb}
\usepackage{color}
\usepackage{ascmac}
\usepackage{url}
\usepackage{ulem}
\usepackage{bm}
\usepackage{bbold}
\usepackage{amsmath}
\usepackage{algorithm}
\usepackage{algpseudocode}
\usepackage{braket}
\usepackage{multirow}
\usepackage{subfiles}
\usepackage{cprotect}
\usepackage{epsfig}
\usepackage{xspace}
\usepackage{tipa}

\newcommand{\HPhi}{ \mathcal{H} \Phi}
\newcommand{\Hphi}{ \mathcal{H} \Phi}

\newcommand{\KOmega}{K$\omega$\xspace}
\newcommand{\Komega}{K$\omega$\xspace}

\newcounter{bla}

\journal{Computer Physics Communications}

\begin{document}

\begin{frontmatter}



\title{\KOmega{} -- Open-source library for the shifted Krylov subspace method of the form $(zI-H)x=b$}


\author[a]{Takeo Hoshi\corref{author}}
\author[b]{Mitsuaki Kawamura}
\author[b]{Kazuyoshi Yoshimi}
\author[b]{Yuichi Motoyama}
\author[b]{Takahiro Misawa}
\author[c]{Youhei Yamaji}
\author[d,b]{Synge Todo}
\author[b]{Naoki Kawashima}
\author[e]{Tomohiro Sogabe}
\cortext[author] {Corresponding author.\\\textit{E-mail address: hoshi@tottori-u.ac.jp} }

\address[a]{Department of Applied Mathematics and Physics, Tottori University,  Tottori-shi, Tottori 680-8552, Japan}
\address[b]{Institute for Solid State Physics, University of Tokyo, Kashiwa-shi, Chiba 277-8581 Japan}
\address[c]{Department of Applied Physics, University of Tokyo, Bunkyo-ku, Tokyo 113-8656, Japan}
\address[d]{Department of Physics, University of Tokyo, Bunkyo-ku, Tokyo  113-0033, Japan}
\address[e]{Department of Applied Physics, Nagoya University, Chikusa-ku, Nagoya 464-8601, Japan}

\begin{abstract}
We develop \KOmega, an open-source  linear algebra library for the shifted Krylov subspace methods. The methods solve a set of shifted linear equations $(z_k I-H)\bm{x}^{(k)}=\bm{b}\, (k=0,1,2,...)$ for a given matrix $H$ and a vector $\bm{b}$, simultaneously.
The leading order of the operational cost is the same as that for a single equation.
The shift invariance of the Krylov subspace
is the mathematical foundation of the shifted Krylov subspace methods.
Applications in materials science are presented to demonstrate the advantages of the algorithm over the standard Krylov subspace methods such as the Lanczos method.
We introduce benchmark calculations of (i) an excited (optical) spectrum and (ii) intermediate eigenvalues by the contour integral on the complex plane.
In combination with the quantum lattice solver $\HPhi$, \KOmega can realize parallel computation of excitation spectra and intermediate eigenvalues for various quantum lattice models.
\end{abstract}


\end{frontmatter}



{\bf PROGRAM SUMMARY}

\begin{small}
\noindent
\KOmega{} -- Open-source library for the shifted Krylov subspace method\\
{\em Authors:} Takeo Hoshi, Mitsuaki Kawamura, Kazuyoshi Yoshimi, Yuichi Motoyama, Takahiro Misawa, Youhei Yamaji, Synge Todo, Naoki Kawashima, Tomohiro Sogabe.
\\
{\em Program Title:} \KOmega [\textipa{k\'ei-\'oumig@}] \\
{\em Journal Reference:}   \\
{\em Catalogue identifier:}                                   \\
{\em Program summary URL:} \\
https://github.com/issp-center-dev/Komega\\
{\em Licensing provisions:} GNU LESSER GENERAL PUBLIC LICENSE Version 3.\\
{\em Programming language:} Fortran~90                                   \\
{\em Computer:} Any architecture with suitable compilers including PCs and clusters.\\
{\em Operating system:} Unix, Linux, macOS.  \\
{\em RAM:} Variable, depending on the dimension of the matrix. \\
{\em Number of processors used:} Arbitrary. 
\\
{\em Classification:4.8 Linear Equations and Matrices. }
\\
{\em External routines/libraries:} BLAS library,
LAPACK library (Used in the sample program),
MPI library (Optional).\\
{\em Nature of problem:
Efficient algorithms, called shifted Krylov subspace algorithms,  designed to  solve the shifted  linear  equations.}\\
{\em Solution method:}
Shifted conjugate gradient method, Shifted conjugate orthogonal conjugate gradient method, shifted bi-conjugate gradient method. \\
{\em Additional comments:} The present paper is accompanied by a frozen copy of \Komega release 2.0.0 that is made publicly available on GitHub (repository https://github.com/issp-center-dev/Komega, commit hash fd5455328b102ec4fa13432496e41c404a0f5a9d).
\\
\\

\end{small}

\section{Introduction}
\label{Intro}

The response of physical observables to an external field or a perturbation is an essential probe in the experimental and theoretical studies of quantum many-body systems.
Theoretically, such a quantity can be formulated as a selected element of the Green's function,
\begin{equation}
\bm{a}^\dagger G(z) \bm{b},
 \label{EQ-GREEN-ELEM}
\end{equation}
with given states (or vectors) $\bm{a}$ and $\bm{b}$, where
$z$ is the complex energy parameter and the Green's function $G(z)$ is defined as the inverse of the shifted Hamiltonian $H$,
\begin{equation}
 G(z)  \equiv (z I - H)^{-1}.
 \label{EQ-GREEN-FN}
\end{equation}
In computational quantum physics, the Hamiltonian $H$ is usually represented by a real-symmetric or Hermitian $M \times M$ sparse matrix.
A typical application of the Green's function formalism is calculating spectra in which the elements $\bm{a}^\dagger G(z_k) \bm{b}$ are calculated at sampling points $z_k$ located near the real axis ($z_k \equiv \omega_k + \mathrm{i} \delta $) with a tiny imaginary part $\delta$
or where the elements $\bm{a}^\dagger G(z_k) \bm{b}$ ($\bm{a} \ne \bm{b}$) are calculated at sampling points $z_k$ for the numerical contour integral on the complex plane.
The dimension $M$ is large, e.g., $M \ge 10^9$, and efficient numerical methods are essential to treat such large matrices.

Formally, the Green's function~(\ref{EQ-GREEN-FN}) can be decomposed as
\begin{equation}
  G(z) = \sum_j \frac{\bm{y}_j \bm{y}_j^\dagger}{z - \lambda_j}
\end{equation}
where $\lambda_j$ and $\bm{y}_j$ are the $j$-th eigenvalue and the corresponding eigenvector of the matrix $H$, respectively, i.e.,
\begin{equation}
H \bm{y}_j = \lambda_j \bm{y}_j.
\label{EQ-ORG-HERM-EV}
\end{equation}
In practice, however, the numerical evaluation of all the eigenpairs is too expensive [the operational cost is $O(M^3)$ and the memory cost is $O(M^2)$] to treat large-scale matrices ($M \ge 10^9$).
In computational quantum physics, a Lanczos-based algorithm combined with the continued fraction technique has been proposed and widely used to calculate the Green's function~\cite{HAYDOCK1980,Dagotto1994}. 
The most expensive operation of the Lanczos-based algorithm is the matrix--vector product, with an operational cost of $O(M^2)$.
Moreover, the sparseness of the Hamiltonian matrix reduces this cost to $O(m)$ where $m$ is the number of nonzero elements in the matrix.
However, the calculation of the Green's function is restricted to diagonal components and it is not straightforward to evaluate the convergences of the obtained results.

In the present paper, we describe the numerical library \KOmega (\url{https://github.com/issp-center-dev/Komega}), which
solves linear equations defined by
\begin{equation}
 (z_k I - H )\bm{x}^{(k)} = \bm{b} \quad (k=0,1,2,\dots,N_\mathrm{eq}-1), 
 \label{EQ-SHIFT-EQ}
\end{equation}
instead of Eq.~(\ref{EQ-ORG-HERM-EV}).
Here, $\bm{b}$ is a given vector, $\{ z_k \}$ are given scalars, $I$ is the identity matrix, and $N_\mathrm{eq}$ is the number of linear equations.  
Eq.~(\ref{EQ-SHIFT-EQ}) is called the shifted linear equation since the matrix on the left hand side is different only by the \lq shift' term, $z_k I$.
In general, the scalar $z$ is a complex number and the matrix of $(z I - H)$ is non-Hermitian.
The solution vector $\bm{x}^{(k)}$ is written as
\begin{equation}
\bm{x}^{(k)} = G(z_k) \bm{b},
 \label{EQ-SOLUTION-VECTOR}
\end{equation}
and satisfies $\bm{a}^\dagger G(z_k) \bm{b} = \bm{a}^\dagger \bm{x}^{(k)}$.
An iterative algorithm is adopted in \KOmega, based on the property that the result of the $n$-th iteration, regardless of $z_k$, can be spanned in the Krylov subspace.
The operational cost is usually dominated by the matrix--vector multiplication procedure for a vector $\bm{v}$ ($\bm{v} \Rightarrow H\bm{v}$).

The library \KOmega is based on efficient algorithms, called shifted Krylov subspace algorithms, designed to efficiently solve the shifted linear equations Eq.~(\ref{EQ-SHIFT-EQ}) with $N_\text{eq}>1$.
The shifted Krylov subspace algorithms and related techniques are briefly explained in this paper 
and the details can be found in~\cite{FROMMER2003, TAKAYAMA2006, YAMAMOTO2007, SOGABE2008, YAMAMOTO2008, TENG2011, Du2011_bIDR, Sogabe2011_EAJAM, TSOGABE2012, DU2012_IDR, Imakura2013_GMRES, GU2014_BOCR, Saito2016, Soodhalte2016, Sun2018} 
and references therein. 
Currently, shifted Krylov subspace algorithms are used in many computational science fields 
such as quantum chromodynamics~\cite{FROMMER2003}, electronic structure calculations~\cite{TAKAYAMA2006,Futamura2010,TENG2011,Takai2017},
excited electron calculations~\cite{GIUSTINO2010}, 
nuclear physics~\cite{MIZUSAKI2010,JIN2017}, 
transport calculations with non-equilibrium Green's function theory~\cite{Iwase2015}, and
nano-structured superconducting systems~\cite{NAGAI2017}.
However, no open-source numerical library/solver of the shifted Krylov subspace methods has yet been developed to our knowledge. 

In this paper, we first review the algorithm implemented in \KOmega briefly in section \ref{Sec:algorithm}. 
Next, the basic information such as installation and usage is introduced in section~\ref{Sec:Usage}. 
Then, some examples using \KOmega as a library or software are illustrated in section~\ref{application}.
Finally, we summarize the paper in section~\ref{sec:summary}.

\section{Algorithms}\label{Sec:algorithm}

\subsection{Overview}
In this section, we briefly review the shifted Krylov subspace methods 
implemented in \KOmega. 
The shifted linear equations of Eq.~(\ref{EQ-SHIFT-EQ})
can be rewritten as
\begin{equation}
 (A + \sigma I)\bm{x}^{\sigma} = \bm{b} 
 \label{EQ-SHIFT-EQ2}
\end{equation}
with $A =  z_0 I - H$ and
$\sigma \equiv z_k - z_0$ for $k=0, 1, 2, \dots, N_\mathrm{eq}-1$. 
Here, the suffix $k$ in Eq.~(\ref{EQ-SHIFT-EQ})
is dropped for simplicity.
Hereafter, Eq.~(\ref{EQ-SHIFT-EQ2}) with $\sigma=0$ 
is called the {\it seed equation},
while the other equations with $\sigma \ne 0$ are called the {\it shifted equations}. 
The accuracy of the approximate solution vector $\bm{x}_n^\sigma$ can be checked 
by monitoring the residual vector
\begin{equation}
 \bm{r}_n^{\sigma} \equiv \bm{b} -  \left( A + \sigma I \right) 
\bm{x}_n^{\sigma},
\label{EQ-RED-VEC}
\end{equation}
where $n$ is the number of iteration steps.

Table~\ref{table:method} shows 
the four solver methods available in \KOmega. 
The methods are classified according to the type
of matrix $A + \sigma I$. 
Users need to select an appropriate method depending on
whether
the matrix $H$ is real-symmetric or Hermitian
and whether the shift constant values $\{z_k\}$ are complex or real\cite{CommentMethod}.

\begin{table*}[tb]
\begin{center}
  \caption{Classification of the Krylov subspace methods.}
  \begin{tabular}{|l|l|l|} \hline
Type of $A + \sigma I$   & Method   & Solver  name\\ \hline \hline
    Real symmetric & Shifted conjugate gradient (CG) method (using real vector) & \verb|komega_cg_r| \\ \hline
    Hermitian & Shifted conjugate gradient (CG) method (using complex vector) & \verb|komega_cg_c|\\ \hline
    Complex symmetric & Shifted conjugate orthogonal conjugate gradient (COCG) method & \verb|komega_cocg| \\ \hline
    Otherwise & Shifted bi-conjugate gradient (BiCG) method &      
\verb|komega_bicg| \\ \hline
  \end{tabular}
  \label{table:method}
  \end{center}
\end{table*}

The mathematical foundation of \KOmega is an iterative algorithm and the numerical solution of Eq.~(\ref{EQ-SHIFT-EQ2}) at the $n$-th iteration is obtained in the Krylov subspace defined as
\begin{equation}
K_n(A, \bm{b}) = \mathrm{span}[\bm{b}, A\bm{b}, A^{2}\bm{b}, \dots, A^{n-2}\bm{b}, A^{n-1}\bm{b}]. 
 \label{EQ-KRYLOV-SUBSPACE2}
\end{equation}
The common mathematical foundation of the shifted Krylov subspace methods is the shift invariance property 
of the Krylov subspace
\begin{equation}
K_n(A + \sigma I, \bm{b}) = K_n(A, \bm{b}),  
 \label{EQ-SHIFT-INVARIANCE}
\end{equation}
and the collinear residual theorem~\cite{FROMMER2003},
which is explained later in this section.

This section explains  
the shifted conjugate gradient (CG) method 
in Table~\ref{table:method}, as an example. 
The other methods,
the shifted conjugate orthogonal CG method and the shifted BiCG method are explained in \ref{appendix-COCG} and \ref{appendix-BiCG}, respectively.

\subsection{Shifted conjugate gradient method}
The shifted CG method is based 
on the CG method~\cite{HestenesStiefel1952} 
that is used for solving $A\bm{x}=\bm{b}$,
when the matrix $A$ is real-symmetric (or Hermitian) and positive definite.

\subsubsection{Seed equation}
The seed equation [Eq.~(\ref{EQ-SHIFT-EQ2}) with $\sigma =0$] is denoted 
as $A\bm{x}=\bm{b}$ and the solution vector at the $n$-th iterative step is denoted as $\bm{x}_{n}$, 
where $\bm{x}_{0}$ is set to be $\bm{x}_{0} = \bm{0}$.
The residual vector at the $n$-th step is denoted as $\bm{r}_{n} = \bm{b}-A\bm{x}_n$. 
The three vectors, $\bm{x}_n$, $\bm{r}_n$, and $\bm{p}_n$, are computed by the following recurrence equations:
\begin{align}
\rho_n &= \bm{r}_n^\dagger \bm{r}_n \label{EQ-REC-RHO} \\
\alpha_n &= \frac{\rho_n}{\bm{p}_n^\dagger A\bm{p}_n} \label{EQ-CG-ALPHA}, \\
\bm{x}_{n+1} &= \bm{x}_{n} + \alpha_{n}\bm{p}_{n}, \\
\bm{r}_{n+1} &= \bm{r}_{n} - \alpha_{n} A \bm{p}_{n} \label{EQ-CG-RN-2REC}, \\
\beta_n &= \frac{\rho_{n+1}}{\rho_{n}}, \label{EQ-REC-BETA}\\
\bm{p}_{n+1} &= \bm{r}_{n+1} + \beta_{n} \bm{p}_{n} \label{EQ-CG-PN-REC}
\end{align}
with the initialization $\bm{p}_0 = \bm{r}_0 = \bm{b}$.

\KOmega uses the three-term recurrence formula~\cite{YAMAMOTO2008} 
\begin{equation}
    \bm{r}_{n+1} = \left(1+\frac{\alpha_n\beta_{n-1}}{\alpha_{n-1}} - \alpha_n A\right) \bm{r}_n - \frac{\alpha_n \beta_{n-1}}{\alpha_{n-1}}\bm{r}_{n-1},
    \label{EQ-CG-RN-3REC}
\end{equation}
which is obtained by
eliminating $\bm{p}_n$ from Eqs.~(\ref{EQ-CG-RN-2REC}) and (\ref{EQ-CG-PN-REC}).
By taking the inner product between $\bm{r}_n$ and Eq.~(\ref{EQ-CG-RN-3REC}), we also obtain
\begin{equation}
    \alpha_n = \frac{\rho_n}{\bm{r}_n^\dagger A\bm{r}_n - \frac{\beta_{n-1}}{\alpha_{n-1}}\rho_n} \label{EQ-REC-ALPHA-N2}.
\end{equation}
These two recurrence equations start with $\beta_{-1}/\alpha_{-1} = 0$ because $\alpha_0 = \rho_0/(\bm{p}_0^\dagger A\bm{p}_0) = \rho_0/(\bm{r}_0^\dagger A\bm{r}_0)$.
The reason for this implementation style will 
be presented in the following subsubsections.

It is crucial that
the $n$-th residual $\bm{r}_n$ belongs to the $(n+1)$-th Krylov subspace $K_{n+1}(A,\bm{b})$ and is orthogonal to the $n$-th Krylov subspace $K_{n}(A,\bm{b})$,
\begin{equation}
    \bm{r}_n \in K_{n+1}(A,\bm{b}) \quad \text{and} \quad \bm{r}_n \perp K_{n}(A, \bm{b}).
\end{equation}
These mean that $\bm{r}_n$ belongs to the one-dimensional orthogonal complementary space of $K_{n}(A,\bm{b})$ within $K_{n+1}(A,\bm{b})$, say $K_{n}^\perp(A,\bm{b})$.

\subsubsection{Shifted equation and collinear residuals}
\label{SEC-CG-SHIFTED}

\begin{figure}
    \centering
    \includegraphics[width=0.6\linewidth]{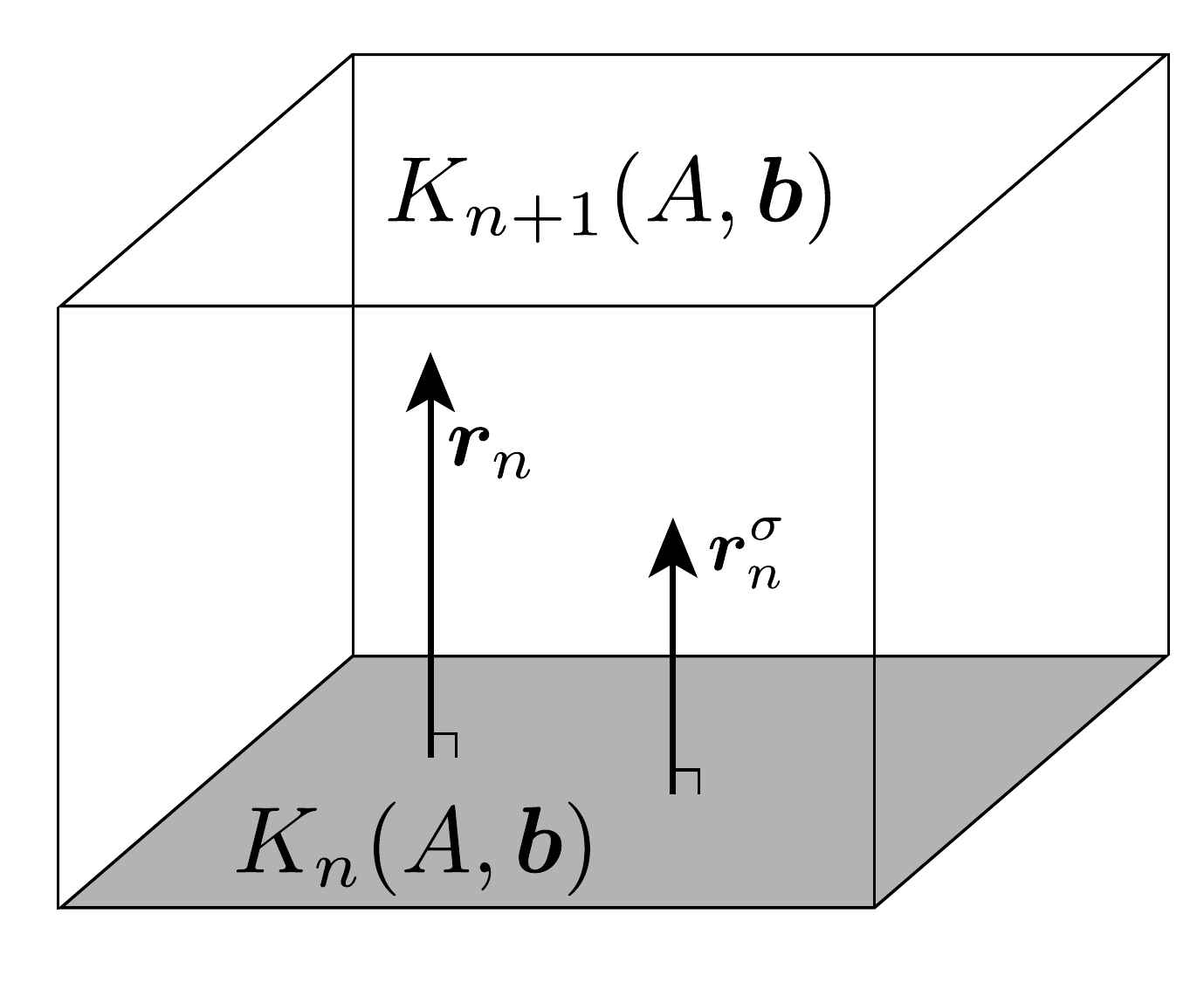}
    \caption{Collinearity of residual vectors. The shaded bottom surface and the cube represent the $n$-th Krylov subspace $K_n(A, \bm{b})$ and the $(n+1)$-th Krylov subspace $K_{n+1}(A, \bm{b})$, respectively.
    The two arrows depict the $n$-th residuals, $\bm{r}_n$ and $\bm{r}_n^\sigma$.
    }
    \label{fig:collinearity}
\end{figure}

When the CG method is applied to a shifted equation 
($\left(A+\sigma I\right)\bm{x} = \bm{b}, \sigma \ne 0$)
with $\bm{x}_0^\sigma = \bm{0}$ and $\bm{r}_0 = \bm{b}$, 
the residual vector $\bm{r}_n^\sigma$ belongs to 
the same one-dimensional subspace $K_n^\perp(A, \bm{b})$ because of Eq.~(\ref{EQ-SHIFT-INVARIANCE}).
Consequently, the residual $\bm{r}_n^\sigma$ 
is parallel to that of the seed equation ($\bm{r}_n$)
\begin{equation}
    \bm{r}_n^\sigma = \frac{1}{\pi_n^\sigma}\bm{r}_n,
    \label{EQ-COLLINEAR}
\end{equation}
where $\pi_n^\sigma$ is a constant~\cite{FROMMER2003} called the collinearity factor.
Figure~\ref{fig:collinearity} is a schematic figure depicting the collinearity ($\bm{r}_n^\sigma \parallel \bm{r}_n)$.

The collinearity factor $\pi_n^\sigma$ 
satisfies the recurrence equation
\begin{equation}
   \pi_{n+1}^\sigma = \left(1+\frac{\alpha_n \beta_{n-1}}{\alpha_{n-1}} + \alpha_n \sigma \right)\pi_n^\sigma - \frac{\alpha_n \beta_{n-1}}{\alpha_{n-1}}\pi_{n-1}^\sigma
\label{EQ-SHIFTEDCG-PI}
\end{equation}
with the initialization $\pi_0^\sigma = \pi_{-1}^\sigma = 1$,
because of Eq.~(\ref{EQ-CG-RN-3REC}) and Eq.~(\ref{EQ-COLLINEAR}).
By substituting $\pi_n^\sigma \bm{r}_n^\sigma$ into $\bm{r}_n$ in Eq.~(\ref{EQ-CG-RN-3REC}),
the recurrence equations for the shifted equations are derived as
\begin{align}
    \alpha_n^\sigma &= \frac{\pi_n^\sigma}{\pi_{n+1}^\sigma}\alpha_n, \label{EQ-SHIFTEDCG-alpha}\\
    \beta_n^\sigma &= \left(\frac{\pi_n^\sigma}{\pi_{n+1}^\sigma}\right)^2 \beta_n, \label{EQ-SHIFTEDCG-beta}\\
    \bm{x}_{n+1}^\sigma &= \bm{x}_{n}^\sigma + \alpha_{n}^\sigma \bm{p}_{n}^\sigma,\label{EQ-SHIFTEDCG-x} \\
    \bm{p}_{n+1}^\sigma &= \frac{1}{\pi_{n+1}^\sigma}\bm{r}_{n+1} + \beta_{n}^\sigma\bm{p}_{n}^\sigma \label{EQ-SHIFTEDCG-p}
\end{align}
with $\bm{x}_0^\sigma = \bm{0}$ and $\bm{p}_0^\sigma = \bm{b}$.
It should be noted that the recurrences 
of Eqs.~(\ref{EQ-SHIFTEDCG-PI}), (\ref{EQ-SHIFTEDCG-alpha}),
(\ref{EQ-SHIFTEDCG-beta}), 
(\ref{EQ-SHIFTEDCG-x}), and  
(\ref{EQ-SHIFTEDCG-p}) 
require no expensive matrix--vector multiplication.

\begin{table}[b]
\begin{center}
  \caption{The operational costs per iteration are listed for (i) the conventional Krylov method, (ii) the shifted Krylov subspace method for the whole elements of the solution vectors, and (iii) the shifted Krylov subspace method for the projection of the solution vectors~\cite{TAKAYAMA2006}. The operational costs of the sparse matrix--vector products (SpMV), the scalar--vector products (SV), and the inner products (Inner) are listed separately. }
  \begin{tabular}{|l|l|l|l|}\hline
   Method & SpMV & SV & Inner  \\  \hline  \hline
   Conventional & $M M_\mathrm{NZ} N_\mathrm{eq}$  & $3 M N_\mathrm{eq}$ & $3 M  N_\mathrm{eq}$  \\  \hline
   Shift (whole) & $M M_\mathrm{NZ} $  & $3 M N_\mathrm{eq}$ & $3 M $  \\ \hline
   Shift (projection) & $M M_\mathrm{NZ} $  & \small{$3 M  + 3 M_\mathrm{left}(N_\mathrm{eq}-1)$} & $3 M $  \\ \hline 
  \end{tabular}
  \label{table:operational_cost}
  \end{center}
\end{table}

\subsubsection{Cost and projection \label{SEC-PROJECTION}}

Table~\ref{table:operational_cost} compares the operational costs of the conventional Krylov method and the shifted Krylov subspace methods~\cite{TAKAYAMA2006}. 
The first row of Table~\ref{table:operational_cost} represents the conventional Krylov method, which solves the $N_\mathrm{eq}$ linear equations independently
in each Krylov subspace.  
Here $M$ and $M_\mathrm{NZ}$ are the dimension and the average number of non-zero elements 
per column of matrix $A$, respectively $(M_\mathrm{NZ} \le M)$. 
The second row represents the cost of
the shifted Krylov subspace method 
for the case where all the elements of the solution vectors $\bm{x}^{\sigma}$ are calculated.  
We find that the operational cost
for the sparse matrix--vector product (SpMV) is drastically reduced ($M M_\mathrm{NZ} N_\mathrm{eq} \rightarrow M M_\mathrm{NZ}$),
since the explicit SpMV appears only in the seed equation $(\sigma=0)$. 
The third row represents the cost of
the shifted Krylov subspace method 
for the case where we do not need all the elements of the solution $\bm{x}^\sigma$, but only its projection, $\bm{y}^\sigma = P \bm{x}^\sigma$ with the $M_\mathrm{left} \times M$ 
projection matrix $P$ ($1 \le M_\mathrm{left} \le M$).
In the third case, we can replace the recurrence equations, Eqs.~(\ref{EQ-SHIFTEDCG-x}) 
and (\ref{EQ-SHIFTEDCG-p}), by
\begin{align}
    \bm{y}_{n+1}^\sigma &= \bm{y}_n^\sigma + \alpha_n^\sigma \bm{u}_n \label{EQ-SHIFTEDCG-y} \\
    \bm{u}_{n+1}^\sigma &= \frac{1}{\pi_n^\sigma}P \bm{r}_n + \beta_n^\sigma \bm{u}_n^\sigma
    \label{EQ-SHIFTEDCG-u}
\end{align}
with $\bm{y}_0^\sigma = \bm{0}$ and $\bm{u}^\sigma_0 = P \bm{b}$.
By this replacement, the number of scalar--vector products (SV) is reduced 
from $3 M N_\mathrm{eq}$ to $3 M  + 3 M_\mathrm{left} (N_\mathrm{eq}-1)$.
For example,  
the calculation of an element of the Green's function
by Eq.~(\ref{EQ-GREEN-ELEM})
\begin{equation}
    G_{ab}(z) \equiv \bm{a}^\dagger \left(zI-H\right)^{-1} \bm{b} = \bm{a}^\dagger \bm{x}(z)
\end{equation}
is a case with $M_\mathrm{left}=1$.
In typical applications, 
both $N_\mathrm{eq}$ and $M_\mathrm{left}$
are much smaller than $M$
($N_\mathrm{eq}, M_\mathrm{left}\ll M$), 
and thus the operational costs in the third row in Table~\ref{table:operational_cost}
are reduced to those in the first row with $N_\mathrm{eq} =1$. 
In other words,
the operational cost in the third case is reduced, typically, to that for solving a single linear equation.

\KOmega implements the third case in Table~\ref{table:operational_cost}, i.e., the shifted Krylov subspace method with projection.
For the seed equation, 
the residual vector $\bm{r}_n$ is updated 
by Eq.~(\ref{EQ-CG-RN-3REC}) and 
the coefficients $\rho_n$, $\alpha_n$, and $\beta_n$
are updated by Eqs.~(\ref{EQ-REC-RHO}),
(\ref{EQ-REC-ALPHA-N2}), and (\ref{EQ-REC-BETA}), respectively. 
For the shifted equations, 
the projected solution vector $\bm{y}_n^\sigma$ 
and the projected search direction vector 
$\bm{u}_n^\sigma$ are updated by Eqs.~(\ref{EQ-SHIFTEDCG-y}) and (\ref{EQ-SHIFTEDCG-u}), respectively,
and the coefficients 
$\pi^\sigma_n$, $\alpha^\sigma_n$, and $\beta^\sigma_n$
are updated by 
Eqs.~(\ref{EQ-SHIFTEDCG-PI}), (\ref{EQ-SHIFTEDCG-alpha}), and
(\ref{EQ-SHIFTEDCG-beta}), respectively.
The solution vectors 
$\bm{x}_{n}^{\sigma}$ for the shifted equations
can be obtained by setting $P=I$ or $M_{\rm left} = M$,
but users should accept the additional memory cost, 
as explained below. 

The present implementation of \KOmega offers a great advantage
not only in the operational cost but also in the memory cost,
because 
an $M$-dimensional vector $\bm{v}$
requires a large memory cost 
in typical applications. 
In the present implementation, 
only three $M$-dimensional vectors,
the residual vectors for the seed equation
($\bm{r}_{n+1}$, $\bm{r}_{n}$, and $\bm{r}_{n-1}$), 
are stored in the memory 
with a memory cost of $O(M^1 M_{\rm eq}^0)$. 
The residual vectors for the shifted equations,
such as $\bm{r}_{n}^\sigma$, 
can be obtained by Eq.~(\ref{EQ-COLLINEAR}) 
with a negligible memory cost 
of $O(M^0 M_{\rm eq}^1)$. 
To store 
the other $M$-dimensional vectors for the shifted equations, 
such as $\bm{x}_{n}^{\sigma}$ and $\bm{p}_n^\sigma$, 
the required memory size will be 
$O(M^1 M_{\rm eq}^1)$, which can be huge.

\subsubsection{Seed switching}
A mathematical technique called seed switching~\cite{YAMAMOTO2008} is adopted 
in \KOmega for an efficient convergence,
because 
the convergence speed of the CG method
can be different among the energy points $\{ \sigma_k \}$.
The residual vectors for several shifted equations 
$(\bm{r}_n^\sigma, \sigma \ne 0)$ can sometimes remain large, 
while that for the seed equation is reduced to be
negligible~\cite{TAKAYAMA2006}.
In this case, we can switch 
the seed equation as $A' = A + \sigma_\mathrm{seed} I$,
where $\sigma_\mathrm{seed}$ is the shift that gives the largest residual $|\bm{r}_n^{\sigma_\mathrm{seed}}|$, or in other words, the smallest collinearity factor $|\pi_n^{\sigma_\mathrm{seed}}|$.
The residual vector for the new seed equation is obtained
by Eq.~(\ref{EQ-COLLINEAR}) 
to be $\bm{r}_n^{\sigma_\mathrm{seed}} = (1/\pi_n^{\sigma_\mathrm{seed}})\bm{r}_n$. 
It is noted that the present implementation style
does not require 
the solution and search direction vectors 
for the new seed equation
($\bm{x}_n^{\sigma_\mathrm{seed}},\bm{p}_n^{\sigma_\mathrm{seed}}$).
\KOmega{} always performs the seed switching after an update.

\subsection{Comparison with other methods}
Here, we briefly describe the merits of calculating $G_{ab}(z)$ by the shifted Krylov subspace method compared to the traditional way based on the Lanczos-based algorithm~\cite{HAYDOCK1980,Dagotto1994}.
In the Lanczos-based algorithm, the Krylov subspace of $K_n(A, \bm{b})$ is generated by the Lanczos-type recurrence formula and the \textit{diagonal} component $G_{bb}$ is given by a continued fraction form.  
In contrast, using the shifted Krylov subspace method, both the \textit{diagonal} components $G_{bb}$ and  the \textit{off-diagonal} components $G_{ab}$ can be calculated directly and simultaneously with the same order of operational cost as that of the Lanczos-based method.
In addition, the accuracy of the obtained results can be evaluated by monitoring the residual vector of Eq.~(\ref{EQ-RED-VEC}). 

The characteristics of the shifted Krylov subspace method is more clarified when compared with (i) independent computation by the standard Krylov subspace method at each single shift point (\lq conventional method' in Table~\ref{table:operational_cost})  on a massively parallel supercomputer and (ii) use of sparse-direct algorithms as in the PEXSI solver~\cite{Lin2014}. As an advantage of the present method, it meets the demand for extremely large problems, such as quantum lattice models with matrix size $M=2^{24} \sim 2^{32}$. Such problems are severe not only in the operational cost but also in the memory cost. In the case of $M= 2^{32} \approx 4.3 \times 10^{9}$, for example, the present method uses three double-precision complex vectors in the recurrence relations and thus requires memory of $3 \times 16 {\rm B} \times M \approx 200 {\rm GB}$, while the others require much larger memory cost. As a limitation of the present \Komega  code, on the other hand, the code is not applicable to the generalized shifted equations $(zB-H) \bm{x} = \bm{b}$ with a positive-definite matrix $B$, unlike the two methods. The generalized shifted equations appear in many physical problems like the electronic state calculation with the overlap matrix $B$. Some of the present authors proposed the shifted Krylov-subspace algorithms for the generalized shifted equations~\cite{TENG2011, TSOGABE2012}, and the implementation of the generalized algorithm in \Komega remains as a future issue.

\section{Usage of \KOmega} \label{Sec:Usage}
In this section, the usage of \KOmega is introduced.
We first introduce how to install \KOmega.
\KOmega provides libraries and a standalone program \verb|ShiftK.out|.
In Sec. \ref{Sec:Flow}, the procedures for using \KOmega are schematically shown. 

\subsection{Installation}
The stable version of \KOmega is distributed on the release page with the Lesser General Public License (LGPL) version 3.
To build  \KOmega a Fortran compiler and a BLAS library are required. 
The following is an example of compiling \KOmega:
\begin{verbatim}
$ ./configure --prefix=install_dir
$ make
$ make install
\end{verbatim}
Here, \verb|install_dir| indicates the full path of the directory where the library is installed. 

\subsection{Schematic flow of \KOmega usage} \label{Sec:Flow}
\begin{figure*}[tb!]
  \begin{center}
    \includegraphics[width=14cm]{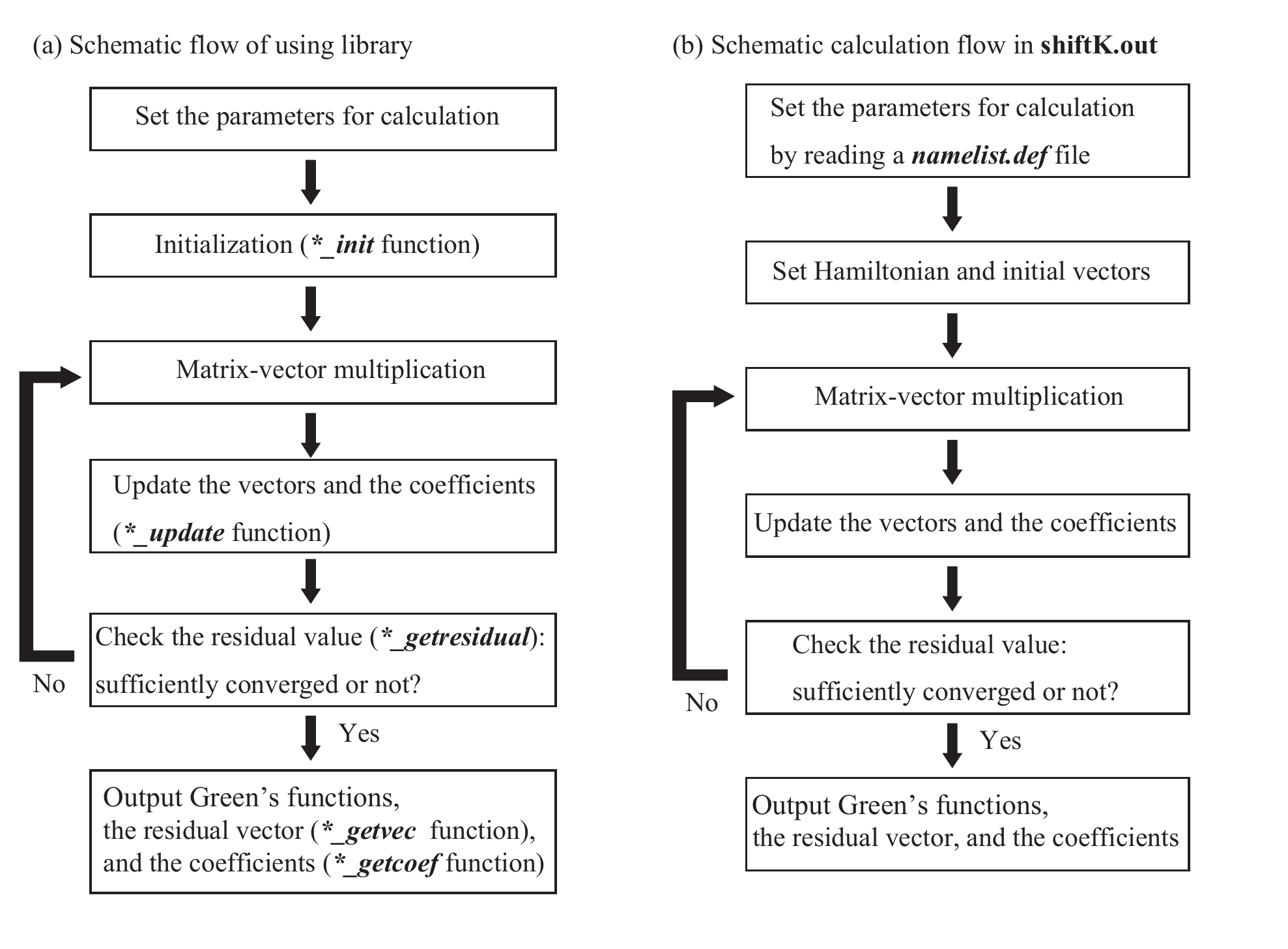}
    \cprotect\caption{(a) Schematic flow of use of library and (b) calculation flow in \verb|shiftK.out|. 
    Here, \verb|*| indicates the name of the solver, such as \verb|komega_cg_r|,  \verb|komega_cg_c|, \verb|komega_cocg|, and \verb|komega_bicg|.}
    \label{fig:flow}
  \end{center}
\end{figure*}

In this subsection, the usage of \KOmega as a library or a standalone program is explained.
In Fig.~\ref{fig:flow}, the schematic flow of the library usage (a) and the corresponding flowchart in the standalone program (b) are shown.

\subsubsection{Library}
\subsubsection* {(i) Preparation of a routine for matrix--vector multiplication \label{RCI}}

\Komega provides a reverse communication interface for 
the matrix--vector multiplication routine ($\bm{v} \Rightarrow H\bm{v}$). The interface requires the preparation of a routine to perform the matrix--vector multiplication, to be called in \Komega.
This interface allows extremely large matrices to be handled,
since the matrix elements are internally generated in the \lq matrix--vector multiplication'  and not stored in the memory. 

\subsubsection* {(ii) Selection of an appropriate solver}
\Komega provides four kinds of numerical solvers.
An appropriate solver should be selected depending on whether the type of Hamiltonian
\({H}\) and the frequency \(z\) are complex or real, as shown in Table~\ref{table:method}.
It is noted that \({ H}\) must be Hermitian or symmetric for a complex or real matrix, respectively.
For an efficient calculation, the seed switching function ~\cite{YAMAMOTO2008} is introduced in all methods.

\subsubsection*{(iii) Calculation}
The calculation is performed in the following steps:\\
(a) Initialization using \verb|*_init| functions to set and initialize internal variables in the library.
For restart calculations, the initial values of the coefficients and the vector should be inputted at this step. \\
(b) Updating of the results iteratively using \verb|*_update| functions, which are called alternately with the matrix--vector product routine in the loop to update the solution. \verb|*_update| also performs the seed switching. At each step, the values of the 2-norm of the residual vector at each shift point are obtained using \verb|*_getresidual| functions. \\
(c) Finalization using \verb|*_finalize| functions to release the memories of the arrays stored in the library.\\
Here, \verb|*| indicates the name of the solver such as \verb|komega_cg_r|,  \verb|komega_cg_c|, \verb|komega_cocg|, and \verb|komega_bicg|.

We summarize this subsection by showing a pseudo code in Algorithm~\ref{ALG-SHIFTEDBICG}, which solves $(zI-H)\bm{x}=\bm{b}$ by the shifted BiCG method (\verb|komega_bicg_*|) implemented in \KOmega.
Users can use the shifted CG method (\verb|komega_cg_*|) or the shifted COCG method (\verb|komega_cocg_*|) in a similar way (they don't require the shadow vectors, $\bm{\tilde{r}}$ and $\bm{\tilde{q}}$.)
The application using the library is shown in Sec.~\ref{SS}.

\begin{algorithm}
\caption{Solve $(zI-H)\bm{x} = \bm{b}$ by Shifted BiCG method with K$\omega$}
\begin{algorithmic}
\Procedure{ShiftedBiCG}{$\bm{x}, H, \{z_i\}, \bm{b}$}
\State ndim $\gets$ dim($H$)
\State nz $\gets$ dim($\{z_i\}$) \Comment{\# of frequencies}
\State x $\gets$ 0 \Comment{initialize output}
\State allocate $r$[ndim] \Comment{residual vector}
\State $r \gets b$
\State allocate $\tilde{r}$[ndim] \Comment{shadow residual vector}
\State $\tilde{r} \gets b^\dagger$
\State allocate $q$[ndim], $\tilde{q}$[ndim] \Comment{storing $Hr$ and $H\tilde{r}$}
\State call \verb|komega_bicg_init|(ndim, nz, ndim, $x$, $\{z_i\}$)
\Repeat
\State $q = Hr$
\State $\tilde{q} = H\tilde{r}$
\State call \verb|komega_bicg_update|($q, r, \tilde{q}, \tilde{r}, x$) \Comment{update and seed switching}
\Until{converged}
\State call \verb|komega_bicg_finilize|()
\State deallocate $r, \tilde{r}, q, \tilde{q}$
\EndProcedure
\end{algorithmic}
\label{ALG-SHIFTEDBICG}
\end{algorithm}

\subsubsection{Standalone program}\label{Sec:Software}
As described in the Introduction, Green's functions are often calculated to investigate the response of physical observables in the quantum many-body systems.
In \KOmega a simple standalone program \verb|ShiftK.out| is provided, which computes the diagonal element of the Green's function:
\begin{equation}
\begin{split}
  G_{aa}(z) = {\bm a}^{\dagger} (z {I} -{H})^{-1}{\bm a}.
  \end{split}\label{dynamicalGreen}
\end{equation}
Here,  ${\bm a}$ is expanded as $\sum_{i} a_i {\bm n}_i$ where ${\bm n}_i$ is the $i$-th basis vector of the Hilbert space. 
An off-diagonal element $G_{ab}$ can be obtained by \verb|ShiftK.out|, when one calculates $G_{aa}, G_{bb}, G_{cc}, G_{dd}$ (${\bm c}\equiv{\bm a}+{\bm b}$, $\bm{d}\equiv{\bm a}+i{\bm b}$), and uses the relation $G_{ab} =\left[(G_{cc}-G_{aa}-G_{bb})+i(G_{dd}-G_{aa}-G_{bb})\right]/2$.\cite{HAYDOCK1980} 
In the following, the usage of \verb|ShiftK.out| is explained.

\subsubsection*{(i) Preparation of an input file}
The input parameters for \verb|ShiftK.out| are categorized into four sections: \verb|cg|, \verb|dyn|, \verb|filename|, and  \verb|ham|.  

The \verb|cg| section sets the numerical condition for the CG (or COCG, BiCG) method.
\verb|maxloops| is the maximum number of iterations and \verb|convfactor| is the threshold value for the convergence criterion of the residual norm, i.e., $\max_\sigma |\bm{r}_n^\sigma| < 10^{-\verb|convfactor|}$.

The \verb|dyn| section specifies the parameters for computing the spectrum.
By setting the parameters \verb|omegamin|, \verb|omegamax|, and \verb|nomega|,
the target frequencies are given as $\omega_i =  \verb|omegamin| + i \times (\verb|omegamax|- \verb|omegamin|)/\verb|nomega|,(i = 0, \dots \verb|nomega|-1)$. 
In \verb|ShiftK.out|, the calculation mode is chosen from \verb|normal|, \verb|recalc|, and \verb|restart| by specifying the parameter \verb|calctype|. \verb|normal| is the mode for computing with the Krylov subspace from scratch. \verb|recalc| is the mode for computing with the Krylov subspace generated in the previous calculation (See (iv) for details).   \verb|restart| is the mode for restarting the calculation from the previous run. 
In the computation of the Green's function, the shifted BiCG or shifted COCG method is automatically selected when $H$ is real or complex, respectively.

\subsubsection*{(a) Input of the matrix of the Hamiltonian}
The \verb|filename| section specifies the name of the input files for the matrix of the Hamiltonian \verb|inham| or the initial excited vector \verb|invec|. The file format of both files is the Matrix Market format~\cite{MatrixMarket}. If \verb|invec| is not specified, a random vector is used as the initial vector. 
An example input file for reading the Hamiltonian matrix and excited vector is as follows:
\begin{verbatim}
&filename
  inham = "Ham.dat"
  invec = "Excited.dat"
/
&cg
  maxloops = 100
  convfactor = 6
/
&dyn
  calctype = "normal"
  nomega = 100
  omegamin = (-2d0, 0.1d0)
  omegamax = ( 1d0, 0.1d0)
/
\end{verbatim}

\subsubsection*{(b) Generation of the matrix of the Hamiltonian}
For the trial use, the matrix of the Hamiltonian is generated in \verb|ShiftK.out| mode. In this mode, the \verb|ham| section is used instead of \verb|inham| in the \verb|filename| section.
In the \verb|ham| section, model parameters are specified to generate the Hamiltonian matrix for a one-dimensional spin chain model:
\begin{equation}
H = \sum_{i=1}^{\verb|nsite|} \left(\begin{array}{ccc}S_{x}^{(i)} & S_{y}^{(i)} & S_{z}^{(i)}\end{array}\right)
\left(\begin{array}{ccc} \verb|Jx| & \verb|Dz| & 0 \\-\verb|Dz| & \verb|Jy| & 0 \\0 & 0 & \verb|Jz|\end{array}\right)
\left(\begin{array}{c}S_{x}^{(i+1)} \\S_{y}^{(i+1)} \\S_{z}^{(i+1)}\end{array}\right),
\end{equation}
where $S^{(i)}_{j}$ is a spin-1/2 operator at a site $i$ with component $j =x, y, z$.
In the \verb|ham| section, the following parameters are specified: the total site of the spin chain \verb|nsite| and the parameters of the strength for spin--spin interactions \verb|Jx|, \verb|Jy|, \verb|Jz|, and \verb|Dz|. An example of using the mode to generate internally the Hamiltonian matrix is as follows:
\begin{verbatim}
&filename
/
&ham
Jx = 1d0
Jy = 1d0
Jz = 1d0
Dz = 1d0
/
&cg
  maxloops = 100
  convfactor = 6
/
&dyn
  calctype = "normal"
  nomega = 100
  outrestart = .TRUE.
/
\end{verbatim}

\subsubsection*{(ii) Run}
After preparing the input file, an executable \verb|ShiftK.out| in terminal is run as follows:
\begin{verbatim}
$ ShiftK.out namelist.def
\end{verbatim}
Here, \verb|namelist.def| is the name of the input file.
The residual values at each step are output to the \verb|residual.dat| file in the working directory.

\subsubsection*{(iii) Results}
After running \verb|ShiftK.out|, the \verb|output| directory is automatically generated. 
In this directory, dynamical Green's functions, the residual vector, and the coefficients are output to \verb|dynamicalG.dat|, \verb|ResVec.dat0|, and \verb|TriDiagComp.dat|, respectively.

\subsubsection*{(iv) Recalculation for additional data (optional)}

The standalone program
\verb|ShiftK.out| provides, as an optional function, 
a recalculation function for additional data,
as explained below. After the successful completion of \verb|ShiftK.out|, users might like to obtain the solution at additional points  $\sigma$ that are not obtained in the completed calculation.
In such cases,
users can calculate the  solution at these points for the shifted equations
with negligible operational costs, 
because the coefficients of $\{ \alpha_n, \beta_n \}$ and the projected residual vectors $\{ P\bm{r}_n \}$ are saved in the files  \verb|TriDiagComp.dat| and \verb|ResVec.dat|, respectively,
and the recurrence relations for the shifted equation Eqs.~(\ref{EQ-SHIFTEDCG-PI}), (\ref{EQ-SHIFTEDCG-alpha}), (\ref{EQ-SHIFTEDCG-beta}), (\ref{EQ-SHIFTEDCG-y}), and (\ref{EQ-SHIFTEDCG-u})
can be solved without any expensive matrix--vector multiplications.

\section{Applications in material science}\label{application}
In this section, we present several numerical results of \KOmega applied to quantum lattice models, to demonstrate typical applications in computational quantum physics. 
The quantum lattice model plays a crucial role for quantum many-body systems and is described by a large sparse matrix $H$. 
The applied studies in this section use the quantum lattice solver package $\HPhi$~\cite{KAWAMURA2017180}. $\HPhi$ was developed as an exact diagonalization solver and can treat a wide range of quantum lattice models, such as the Hubbard model, the Kondo model, and the Heisenberg model. \KOmega can be called from $\HPhi$ and the users of $\HPhi$ can use the shifted Krylov subspace solvers. Here, three typical examples are explained. 

\subsection{Calculation of internal eigenpairs}~\label{SS}
The first example is the calculation of internal eigenvalues by a contour-integral method~\cite{SAKURAI2003}.
The method is an efficient method for obtaining eigenpairs (eigenvalues and eigenvectors) in a specified eigenvalue range.
Although a number of software packages for this method have already been developed (for example, zPares~\cite{zpares} and FEAST~\cite{FEAST,FEAST_PRB}), this example illustrates how to use K$\omega$ in the contour-integral method.

Here, we briefly explain how to obtain eigenpairs by using the  method
following ref.~\cite{SS_JPSJ2013}.
In the contour-integral method, the projection matrix on the $n$-th eigenvector,
\begin{align}
{P}_{n}=\bm{y}_n \bm{y}_n^{\dagger},
\end{align}
plays a key role.
By multiplying ${P}_{n}$ by a vector
$\bm{\phi}=\sum_n a_n \bm{y}_n$, we can extract 
the component of $\bm{y}_n$ as
\begin{align}
{P}_{n}\bm{\phi}=a_{n}\bm{y}_n.
\end{align}
On the other hand, the
projection matrix ${P}_{n}$ can be expressed by an 
integration in the complex plane as
\begin{align}
&{P}_{\Gamma}
= \sum_{\lambda_n\in \Gamma}{P}_{n}
=\frac{1}{2\pi\mathrm{i}}\oint_{\Gamma} \frac{1}{zI-{H}}dz,
\end{align}
where $\Gamma$ represents a contour on the complex plain.

Using $P_{\Gamma}$,
we can extract only the eigenvectors 
whose eigenvalues exist within $\Gamma$ as
\begin{align}
\bm{s}_{0,0} &= {P}_{\Gamma}\bm{\phi_{0}}
=\frac{1}{2\pi\mathrm{i}}\oint_{\Gamma} \frac{1}{zI-{H}}\bm{\phi_{0}}dz  \notag \\ 
&= \frac{1}{2\pi\mathrm{i}}\oint_{\Gamma} \bm{\phi_{0}^{\prime}}dz 
=\sum_{\lambda_n\in \Gamma}a_{n}\bm{y}_n. \label{eq:calc_s00}
\end{align}
In this process,
we can use K$\omega$ to obtain $ \bm{\phi^{\prime}_{0}}=(zI-H)^{-1}\bm{\phi_{0}}$.

The $k$-th
moment can be simply expressed as follows:
\begin{align}
&\bm{s}_{k,0}=({H}-z_{0}I)^{k}{P}_{\Gamma}\bm{\phi_{0}}
=\frac{1}{2\pi\mathrm{i}}\int_{\Gamma} \frac{(z-z_{0})^{k}}{zI-{H}}\bm{\phi_{0}} dz,
\end{align}
where $z_{0}$ is an arbitrary complex number.
Thus, once we obtain $\bm{\phi_{0}^{\prime}}=(zI-H)^{-1}\bm{\phi_{0}}$,
we can calculate the Krylov subspace defined as
\begin{align}
{K}_{N_{k}\times 1}={\rm span}[\bm{s}_{0,0},\bm{s}_{1,0},\dots,\bm{s}_{N_k-1,0} ],
\end{align}
where $N_{k}$ represent the number of moments.
Taking other vectors $\bm{\phi}_{l}$ $(l = 1\cdots N_l -1)$,
we can extend the Krylov subspace as
\begin{align}
{K}_{N_{k}\times N_{l}}
={\rm span}[&\bm{s}_{0,0},\dots,\bm{s}_{N_k-1,0}; \notag \\
&\bm{s}_{0,1},\dots,\bm{s}_{N_k-1,1}; \dots \notag \\
&\bm{s}_{0,N_{l}-1},\dots,\bm{s}_{N_k-1,N_{l}-1} ].
\end{align}
Then, the matrix $S$ with vectors $\bm{s}_{k,l}$ is constructed as
\begin{align}
S=(&\bm{s}_{0,0},\dots,\bm{s}_{N_k-1,0}, \notag \\
&\bm{s}_{0,1},\dots,\bm{s}_{N_k-1,1}, \dots \notag \\
&\bm{s}_{0,N_{l}-1},\dots,\bm{s}_{N_k-1,N_{l}-1}).
\end{align}
By performing singular value decomposition for $S$,
we obtain
\begin{align}
&S=U\Sigma V^{\dagger}.
\end{align}
The number of non-zero singular values in $\Sigma$ is just that of the independent vectors ($M_\mathrm{nz}$) in $N_{N_{k}\times N_{l}}$ 
spanned by $\bm{s}_{k,l}$.
By using only the left-singular vectors with non-zero singular values,
we can construct the matrix $\tilde{U}$,
\begin{align}
\tilde{U}=(\bm{u}_{0},\dots,\bm{u}_{M_\mathrm{nz}-1}).
\end{align}
Using $\tilde{U}$,
we project the original Hamiltonian
on the small matrix (dimension $M_\mathrm{nz}\times M_\mathrm{nz}$) whose eigenvalues 
exist within $\Gamma$ as
\begin{align}
&\tilde{{H}}=\tilde{U}^{\dagger} {H}\tilde{U}.
\end{align}
By diagonalizing $\tilde{H}$,
we can obtain all the eigenvalues and eigenvectors which exist within $\Gamma$.

As an example, 
we apply the contour-integral method
to the $12$-site Heisenberg chain model defined in Eq.~(\ref{eq:heisenberg}).
For simplicity, we obtain the eigenvalues around the ground state (the lowest eigenvalue).
To perform the integration on the complex plane,
we use the following points:
\begin{align}
z_{j}=\gamma+\rho\exp\Big[\frac{2\pi\mathrm{i}}{N_z}\left(j+\frac{1}{2}\right)\Big],
\end{align}
where we take $\gamma=-5$,
$\rho=0.8$, and $N_{z}=100$ as an example.

In Table~\ref{table:SS}, we show the result of 
the method for the $L=12$ Heisenberg chain,
whose Hilbert dimension is $924$.
To see the convergence of the contour-integral method,
we show the obtained eigenvalues for
several different conditions. 
By taking about $50$ different 
bases, we can obtain the correct eigenvalues
including degeneracies.
To perform the above calculations, the matrix data file in the Matrix Market format (\verb|Ham.dat|) is generated using $\HPhi$ as in \ref{Sec:HPhi}. Using \verb|Ham.dat|, we can perform the contour-integral method by K$\omega$. Two sample codes are provided in the web site\cite{SampleCodeSS}: One is a simple code for obtaining ${\bm \phi}_0'$ in Eq.~(\ref{eq:calc_s00}) by BiCG method (see Algorithm 1) and the other is for performing the contour-integral method.
\begin{table}[tb!]
\caption{Eigenvalues obtained by $\mathcal{H}\Phi$ 
(full diagonalization using LAPACK) 
and the contour-integral method. 
For the contour-integral method, we show results 
for several different $N_k$ and $N_l$,
denoted as SS($N_{k}$,$N_{l}$).
When the rank of SS($N_{k}$,$N_{l}$) is smaller than $n$,
we do not obtain the $n$-th eigenvalues.
This case is indicated as NA (not available) in the table.
In the actual calculation, we discard singular values that are smaller 
than $10^{-3}$.}
\begin{center}
\begin{tabular}{llllll}
\hline 
$E_{n}$ & $\mathcal{H}\Phi$   &  SS(10,1)  & SS(10,2)   & SS(10,5)   \\ \hline
0       & $-5.387391$           &  $-5.387391$  & $-5.387391$  & $-5.387391$  \\
1       & $-5.031543$           &  $-5.031543$  & $-5.031543$  & $-5.031543$  \\
2       & $-4.777389$           &  $-4.777389$  & $-4.777389$  & $-4.777389$  \\
3       & $-4.569374$           &  $-4.569374$ & $-4.569374$  & $-4.569374$  \\
4       & $-4.569374$           &  $-4.297689$        & $-4.569374$  & $-4.569374$  \\
5       & $-4.297689$           &  NA        & $-4.297689$  & $-4.297689$  \\
6       & $-4.297689$           &  NA        & $-4.297689$  & $-4.297689$  \\
\hline
\end{tabular}
\end{center}
\label{table:SS}
\end{table}

\subsection{Calculation of spectra using ShiftK.out}\label{dynamical}

The second example is a spectral calculation of Eq.~(\ref{dynamicalGreen}) by the standalone program \verb|ShiftK.out| explained in Sec.~\ref{Sec:Software}.
A typical application for a quantum lattice model is the calculation of 
the dynamical correlation factor, since it is often used to analyze low-energy structures. 
The dynamical correlation factor is defined as
\begin{align}
D(\omega,\eta)=-\frac{1}{\pi}\mathrm{Im}
\Big[
\bm {\phi_{0}}^{\dagger}
A
\left[ (\omega-E_{0}-\mathrm{i}\eta)I - H\right]^{-1}
B
\bm {\phi_{0}}
\Big],
\end{align}
where $\bm {\phi_{0}}$ is the ground state vector,
$A$ and $B$ are the matrices to generate excited states, $\omega$ represents the frequency,
and $\eta$ represents the smearing factor.
For example,
by taking $A=B=S^{z}$,
we can calculate the dynamical spin structure factors,
which can be directly measured by the neutron scattering experiments.
The key part in the calculation of the dynamical 
correlation factors is solving the linear equation defined as
\begin{align}
\bm{\phi}^{\prime} 
=\left[(\omega-E_{0}-\mathrm{i}\eta)I - H \right]^{-1} B \bm{\phi}_{0}.
\end{align}

As an example of a calculation of dynamical Green's functions by \verb|ShiftK.out|, 
we present the calculation of the dynamical spin structure factors on a one-dimensional Heisenberg chain, whose Hamiltonian is defined as
\begin{align}
H=J\sum_{i=0}^{L-1}\bm{S}_{i}\cdot \bm{S}_{i+1}, \label{eq:heisenberg}
\end{align}
where $\bm{S}$ represents the
spin 1/2 operator and we take 
a magnetic interaction of $J=1$ and system size of $L=12$.
The matrix data file of the Hamiltonian $H$ is generated by $\HPhi$ in the Matrix Market format.
The dynamical spin structure factors are defined as
\begin{equation}
S(q, \omega,\eta)= - \frac{1}{\pi}\mathrm{Im} \Big[ \bm{\phi}_{0}^\dagger S^z(-q) \left[(\omega-E_{0}-\mathrm{i}\eta)I-H \right]^{-1} S^z(q) {\bm\phi}_{0} \Big], \label{eq:dynamicalspin}
\end{equation}
where $S^z(q)=\sum_{j=0}^{L-1} e^{i q j}S^z_{j} $ and $q$ is the wave vector. The input file for  \verb|ShiftK.out| is given as follows ($E_0=0, -5.5\le \omega \le 0.0$, $ \eta = 0.02$):
\begin{verbatim}
&filename
  inham = "Ham.dat"
  invec = "excitedvec.dat"
/
&cg
  maxloops = 1000
  convfactor = 6
/
&dyn
  calctype = "normal"
  nomega = 1000
  omegamin = (-5.5, -0.02d0)
  omegamax = ( 0.0, -0.02d0)
/
\end{verbatim}
Here, \verb|Ham.dat| and \verb|excitedvec.dat| are the input files for the Hamiltonian matrix $H$ and the excited vector $S^z(q) \bm{\phi}_{0}$. These files can be  obtained using the quantum lattice solver $\HPhi$. 
The details are shown in \ref{Sec:HPhi}. 

\begin{figure}[tb!]
  \begin{center}
    \includegraphics[width=0.8\columnwidth]{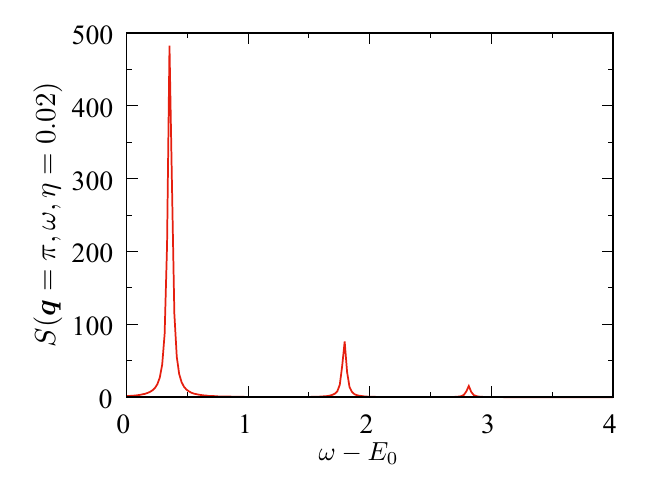}
    \caption{Dynamical Green's function for the one-dimensional Heisenberg model.
    In this calculation, we take $\eta=0.02$, $\omega_\mathrm{min}=-5.5$, $\omega_\mathrm{max}=0.0$,
    and the energy shift $E_{0}=-5.387$, which is the energy of the ground state.}
    \label{fig:szqq}
  \end{center}
\end{figure}

In Fig.~\ref{fig:szqq}, we show the numerical result of $S( q = \pi, \omega, \eta)$. Here, we shift $\omega$ by the ground state energy $E_{0} = -5.387$.
The frequencies where $S(q = \pi, \omega, \eta)$ have sharp peaks correspond to the excitation energies induced by $S_{z}(q=\pi)$. The convergence of the residuals can be checked in the standard output. 

\subsection{Calculation of optical conductivity using $\HPhi$}~\label{optical}

The third example is the calculation of a spectrum or dynamical correlation function using $\HPhi$, demonstrating that using K$\omega$ in $\HPhi$ can be a powerful tool for various problems. 
As an example, we show the optical conductivity in the extended Hubbard model, which is a typical model for strongly correlated electron systems.
We note that the optical conductivity is often used for 
examining the metallic or insulating behavior of 
correlated electron systems.
The optical conductivity can be calculated from
the current--current correlation 
$I(\omega,\eta)$, which is defined as
\begin{align}
&j_{x}=\mathrm{i}\sum_{i,\sigma}
(
c_{\bm{r}_{i}+\bm{e}_{x},\sigma}^{\dagger}c_{\bm{r}_{i},\sigma}
-c_{\bm{r}_{i},\sigma}^{\dagger}c_{\bm{r}_{i}+\bm{e}_{x},\sigma}
), \\
&I(\omega,\eta)=\mathrm{Im}
\Big[
\bm{\phi_{0}}^{\dagger}
j_{x}
[H-(\omega-E_{0}-\mathrm{i}\eta)I]^{-1}
j_{x}
\bm{\phi_{0}}
\Big],
\end{align}
where $\bm{e}_{x}$ is the unit translational vector
in the $x$ direction.
The ground state vector $\bm{\phi_{0}}$ is calculated by an exact diagonalization solver built in $\HPhi$.
To obtain $\bm{\phi_{0}}^{\prime}=[H-(\omega-E_{0}-\mathrm{i}\eta)I]^{-1}j_{x}\bm{\phi_{0}}$, we use K$\omega$.
From the current--current correlation,
the regular part (i.e., without the Drude part at $\omega=0$) 
of the optical conductivity 
is defined as
\begin{align}
\sigma_\mathrm{reg}(\omega)=\frac{I(\omega,\eta)+I(-\omega,-\eta)}{\omega N_{s}},
\end{align}
where $N_\mathrm{s}$ is the number of sites.

To directly compare the optical conductivity 
with previous studies~\cite{MerinoPRB2005}, 
we calculate the optical conductivity in the 
extended Hubbard model defined as
\begin{align}
H=&-t\sum_{\langle i,j\rangle}
\Big[c_{i\sigma}^{\dagger}c_{j\sigma}+
c_{j\sigma}^{\dagger}c_{i\sigma}\Big]+
U\sum_{i}n_{i\uparrow}n_{i\downarrow} \notag \\
&+ V\sum_{\langle i,j\rangle}N_{i}N_{j}+
V^{\prime}\sum_{\langle\langle i,j\rangle\rangle}N_{i}N_{j},
\end{align}
where $c_{i\sigma}$ and $c^{\dagger}_{i\sigma}$ denote the annihilation and creation operators of an electron
at site $i$ with spin $\sigma$, and $n_{i\sigma}=c^{\dagger}_{i\sigma}c_{i\sigma}$ represents the number operator of an electron at
site $i$ with spin $\sigma$.
We note that this model at quarter filling is 
an effective model
for organic conductors.
We perform a calculation for an $N_\mathrm{s}=4\times4$ system with
 $t=1$, $U=4$, $V=3$, and $V^{\prime}=5$. 

\begin{figure}[tb!]
  \begin{center}
    \includegraphics[width=0.7\columnwidth]{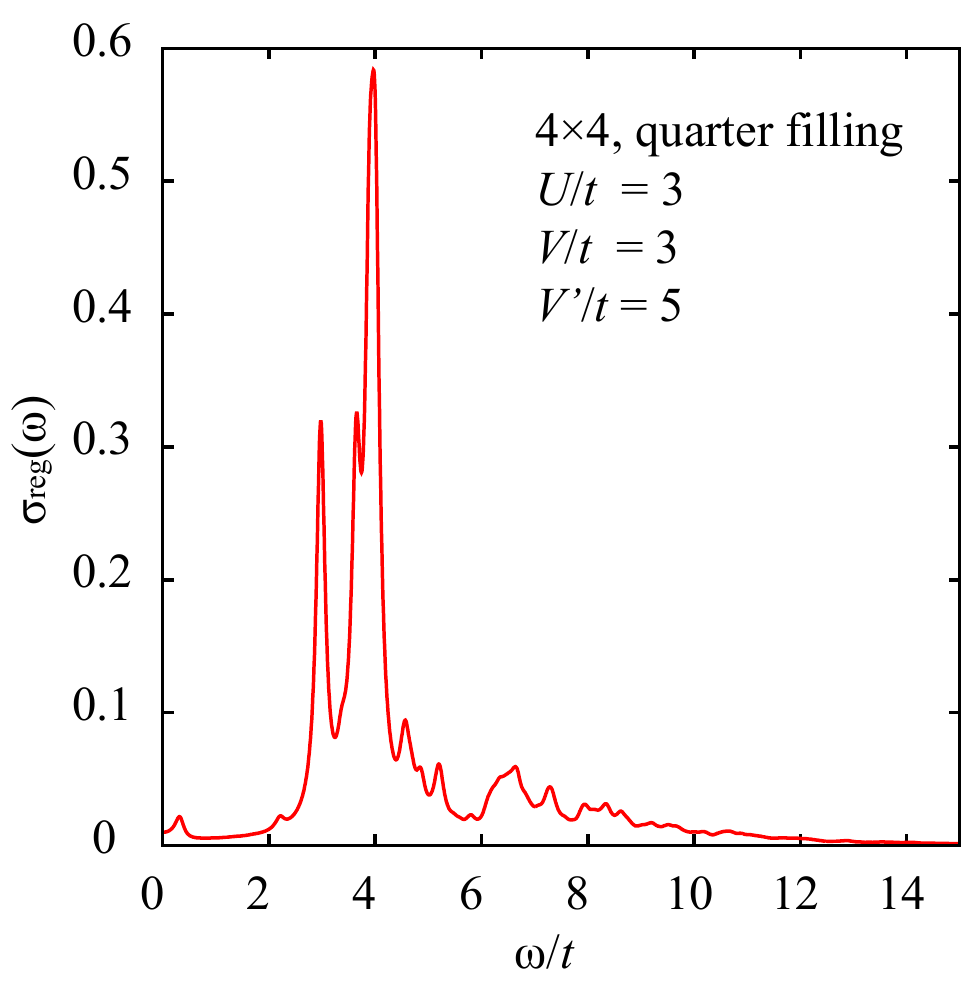}
    \caption{Optical conductivity for the extended Hubbard model.
    In this calculation, we take $\eta=0.1$, $\omega_\mathrm{min}=0$, $\omega_\mathrm{max}=15.0$,
    and the number of $\omega$ as $2000$.}
    \label{fig:opt}
  \end{center}
\end{figure}

\begin{figure}[tb!]
  \begin{center}
    \includegraphics[width=0.7\columnwidth]{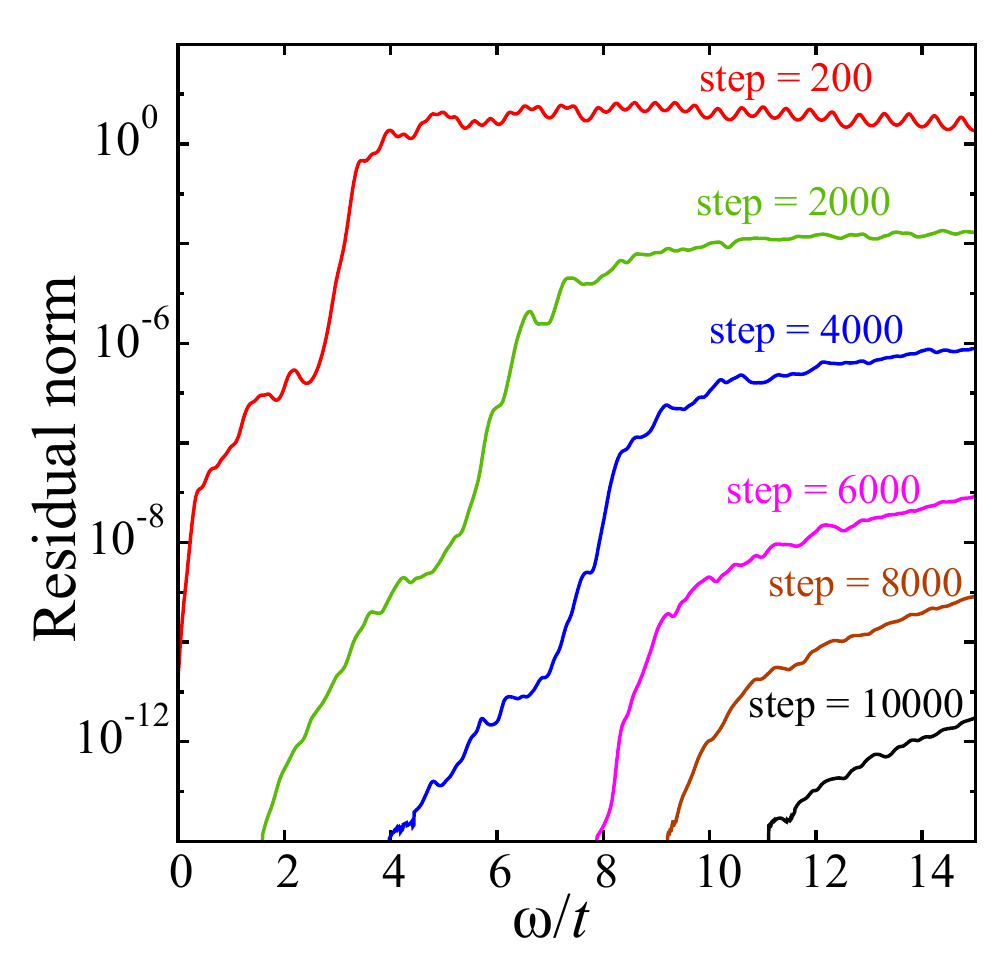}
    \caption{Residual norm $\max_{\sigma}|{\bm r}_n^{\sigma}|$ of $I(\omega,\eta)$ for several iteration steps.
    For the extended Hubbard model, we find that
    the residual for the high-energy part remains large.
     This result indicates that it is necessary to perform more than $10^4$ iterations to obtain a well-converged correlation function when the threshold value is set to $10^{-12}$. 
    }
    \label{fig:res}
  \end{center}
\end{figure}

In Fig.~\ref{fig:opt},
we show the result of the 
optical conductivity in the extended Hubbard model.
This result is consistent 
with a previous study~\cite{MerinoPRB2005}.
We also show the residuals
of $I(\omega,\eta)$ in Fig.~\ref{fig:res}.
By calculating the residuals,
we can check whether 
the obtained dynamical correlation functions
are well converged.
As mentioned above, 
this is one of the main advantages of the shifted Krylov subspace method.
Furthermore,
by examining the $\omega$ dependence of the
residuals we can also identify 
bottlenecks for the convergence.

It is noteworthy that the matrix size becomes huge for the calculation of realistic models. In such cases, the matrix--vector product dominates the numerical cost of the shifted-Krylov method. To address this problem, we designed the K$\omega$ library to not directly include the matrix--vector product in the function, as explained in Sec.~\ref{RCI}. If users can prepare the matrix--vector product function with high parallelization efficiency, the calculation time is greatly suppressed. In fact, the calculation of dynamical Green's functions by $\Hphi$ shows high parallelization efficiency, since the matrix--vector product is well parallelized~\cite{KAWAMURA2017180}. From this point of view, this library is considered to be suitable for large-scale calculations. 

\section{Summary}
\label{sec:summary}
We developed the numerical library \Komega
for solving shifted linear equations with the shifted Krylov subspace methods. 
The present paper details applications for quantum many-body problems that appear in condensed matter physics. As a demonstration, numerical results including dynamical Green's functions, optical conductivity, and eigenvalues through the contour-integral method are shown. The method is also applicable to other computational physics areas, and hence \Komega could be a useful numerical library in computational physics. 

\section{Acknowledgements}
TH wishes to thank S. Yamamoto for helpful advice regarding the shifted Krylov subspace method.
KY and TM wish to thank M. Naka and H. Seo for fruitful discussion about calculation of the optical conductivity.
This work was supported in part by the Ministry of Education, Culture, Sports, Science and Technology (MEXT) of Japan as a social and scientific priority issue (Creation of new functional devices and high-performance materials to support next-generation industries; CDMSI) to be tackled by using post-K computer.
This work was also supported in part by the KAKENHI funds (20H00581, 19H04125, 17H02828).
We would also like to acknowledge support from 
the ``\textit{Project for advancement of software usability in materials science (PASUMS)}" 
by the Institute for Solid State Physics, University of Tokyo, 
for the development of \Komega ver.\ 1.0.
This work was also supported by the Building of Consortia for the Development of Human Resources in Science and Technology from the MEXT of Japan.
We also wish to thank the Supercomputer Center of Institute for Solid State Physics, University of Tokyo for the use of their numerical resources.

\appendix

\section{Shifted conjugate orthogonal conjugate gradient method}\label{appendix-COCG}
When all the shifted matrices $A + \sigma I$ are complex symmetric, users should use the shifted conjugate orthogonal conjugate gradient (COCG) method~\cite{YAMAMOTO2008}.
In the COCG method~\cite{vanderVorst1990}, 
recurrence equations can be obtained by changing the inner products from $\bm{x}^\dagger\bm{y}$ to  $\bm{x}^\mathrm{T}\bm{y}$ in those for the CG method and in the definition of the orthogonality and $A$ conjugate.
The $n$-th residual $\bm{r}_n$ is orthogonal (defined by the ordinary inner product) not to $K_{n}(A,\bm{b})$ but to $K_{n}(A,\bm{b})^*$.
The collinearity of the shifted residuals still holds~\cite{YAMAMOTO2008}, however, since $\bm{r}_n$ belongs to the one-dimensional orthogonal complementary space of $K_{n}(A,\bm{b})^*$ within $K_{n+1}(A,\bm{b})$.
In the result, the shifted equations have the same form as those for the shifted CG method.

\section{Shifted bi-conjugate gradient method}\label{appendix-BiCG}
When $A+\sigma I$ is not symmetric or Hermitian, users should use 
the shifted bi-conjugate gradient (BiCG) method~\cite{FROMMER2003}.
\subsection{Bi-conjugate gradient method}
The BiCG method~\cite{Fletcher1976} generates sequences of approximate solutions $\bm{x}_n$, residuals $\bm{r}_n$, shadow residuals $\tilde{\bm{r}}_n$, search directions $\bm{p}_n$, and shadow search directions $\tilde{\bm{p}}_n$ so that
the pair of residuals and shadow residuals forms a bi-orthogonal system, that is,
\begin{equation}
\tilde{\bm{r}}_i^\dagger\bm{r}_j = 0 \quad (i\ne j); \quad \tilde{\bm{r}}_i^\dagger \bm{r}_i \ne 0.
\end{equation}
These vectors are iteratively generated by the following recurrence equations, as for the CG method, 
\begin{align}
\rho_n &= \tilde{\bm{r}}_n^\dagger \bm{r}_n, \\
\alpha_n &= \frac{\rho_n}{\tilde{\bm{p}}_n^\dagger A\bm{p}_n}, \\
\bm{x}_{n+1} &= \bm{x}_{n} + \alpha_{n}\bm{p}_{n}, \\
\bm{r}_{n+1} &= \bm{r}_{n} - \alpha_{n} A \bm{p}_{n}, \label{EQ-RN-2REC} \\
\tilde{\bm{r}}_{n+1} &= \tilde{\bm{r}}_{n} - \alpha_{n}^* A^\dagger \tilde{\bm{p}}_{n}, \\
\beta_n &= \frac{\rho_{n+1}}{\rho_{n}}, \\
\bm{p}_{n+1} &= \bm{r}_{n+1} + \beta_{n}  \bm{p}_{n}, \label{EQ-PN-REC} \\
\tilde{\bm{p}}_{n+1} &= \tilde{\bm{r}}_{n+1} + \beta_{n}^*  \tilde{\bm{p}}_{n},
\end{align}
with the initialization $\bm{r}_0 = \bm{p}_0 = \bm{b} - A\bm{x}_0$ and $\tilde{\bm{r}}_0 = \tilde{\bm{p}}_0 = \tilde{\bm{b}}$, where $\bm{x}_0$ in an input vector and
$\tilde{\bm{b}}$ is an arbitrary vector that is not orthogonal to $\bm{b}$.

The residuals $\bm{r}$ and $\tilde{\bm{r}}$ and the coefficient $\alpha$ follow three term recurrences:
\begin{align}
    \bm{r}_{n+1} &= \left(1+\frac{\alpha_n\beta_{n-1}}{\alpha_{n-1}} - \alpha_n A\right) \bm{r}_n - \frac{\alpha_n \beta_{n-1}}{\alpha_{n-1}}\bm{r}_{n-1}, \\
    \tilde{\bm{r}}_{n+1} &= \left(1+\left(\frac{\alpha_n\beta_{n-1}}{\alpha_{n-1}}\right)^* - \alpha_n^* A^\dagger \right) \tilde{\bm{r}}_n - \left(\frac{\alpha_n \beta_{n-1}}{\alpha_{n-1}}\right)^*\tilde{\bm{r}}_{n-1}, \\
    \alpha_n &= \frac{\rho_n}{\tilde{\bm{r}}_n^\dagger A\bm{r}_n - \frac{\beta_{n-1}}{\alpha_{n-1}}\rho_n},
\end{align}
with $\beta_{-1}/\alpha_{-1} = 0$.
\Komega{} adopts these three term recurrences.

\subsection{Shifted equation}
When starting the BiCG iteration for the seed equation $A\bm{x}=\bm{b}$ and some shifted equation $(A+\sigma I)\bm{x} = \bm{b}$ with the same initial guess $\bm{x}_0 = \bm{x}_0^\sigma = \bm{0}$ and shadow residuals $\tilde{\bm{r}}_0$,  the $n$-th residuals $\bm{r}_n$ and $\bm{r}_n^\sigma$ are parallel to each other~\cite{FROMMER2003},
\begin{equation}
    \bm{r}_n^\sigma = \frac{1}{\pi_n^\sigma} \bm{r}_n.
\end{equation}
The approximate solutions for the shifted equations are calculated via the same equations in the shifted CG method, Eqs.~(\ref{EQ-SHIFTEDCG-PI})-(\ref{EQ-SHIFTEDCG-p}).

\section{Generate input files by $\Hphi$} \label{Sec:HPhi}
This appendix shows the Hamiltonian matrix in the Matrix Market format and the excited vector generated using $\Hphi$. These files can be read by \verb|shiftK.out|. As an example, we focus on the $12$-site one-dimensional Heisenberg chain, whose Hamiltonian is defined in Eq.~(\ref{eq:heisenberg}) and the excited vector $S^z(q) \bm{\phi}_{0}$ in Eq.~(\ref{eq:dynamicalspin})

\subsection{Output Hamiltonian matrix}
An input file for $\mathcal{H}\Phi$ is given as follows:
\begin{verbatim}
L = 12
model = "Spin"
lattice = "chain"
method = "FullDiag"
J = 1.0
2Sz = 0
HamIO = "out"
\end{verbatim}
Using this input file,
we can obtain the output file of the Hamiltonian matrix in the Matrix Market format (\verb|zvo_Ham.dat| renamed as \verb|Ham.dat| in Sec.~\ref{application}).

\subsection{Output the excited vector}
To obtain the excited vector $S^z(q) \bm{\phi}_{0}$ by $\Hphi$,
we first calculate the ground state $\bm{\phi}_{0}$  and
then calculate the excited vector by multiplying $\bm{\phi}_{0}$  by $S^z(q)$.
The ground state  $\bm{\phi}_{0}$ can be obtained by the following input file:
\begin{verbatim}
L  = 12
model = "Spin"
lattice = "chain"
method = "CG"
J = 1.0
2Sz = 0
EigenVecIO = "Out"
\end{verbatim}
The input file for calculating the excited state $S^z(q) {\bm \phi}_{0}$ is given as follows:
\begin{verbatim}
L = 12
model = "Spin"
lattice = "chain"
method = "CG"
J = 1.0
2Sz = 0
CalcSpec="Normal"
SpectrumQL = 0.5
EigenVecIO ="In"
OutputExcitedVec = "Out"
\end{verbatim}
Here, the wave vector is specified using \verb|SpectrumQL|. In the above case, $S^z(q = \pi){\bm \phi}_0$ is obtained. For details, see $\Hphi$ user's manual. After the calculations are finished, the output file of the excited vector \verb|zvo_excitedvec_rank_0.dat| (renamed as \verb|excitedvec.dat| in Sec.~\ref{application}) is obtained.






\bibliographystyle{apsrev}
\bibliography{mybibfile}

\begin{thebibliography}{38}
\bibitem[{Com()}]{CommentMethod}
\bibinfo{note}{It is noted that the method for real-symmetric matrix $H$ can be
  applied to the problem for complex Hermitian matrix $H$~\cite{JIN2017},
  because the eigenvalue problem of Eq.~(\ref{EQ-ORG-HERM-EV}) with an $M
  \times M$ complex Hermitian matrix $H$ is transformed into one with a $2M
  \times 2 M$ real-symmetric matrix.}
\end{thebibliography}







\end{document}